\documentclass[12pt,leqno]{article}
\usepackage[pdftex]{graphicx}
\usepackage{amsfonts, color}

\newcommand{\ind}{\hspace{.25in}}

\newcounter{conjecture}
\newcounter{remark}
\newcounter{corollary}
\newenvironment{conjecture}{\medskip{\bf Conjecture.\ }\em}{\rm}

\newtheorem{corollary}{Corollary}
\newtheorem{theorem}{Theorem}
\newtheorem{lemma}{Lemma}
\newtheorem{proposition}{Proposition}

\newcommand{\aaa}{\alpha}

\newcommand{\lll}{\label}
\newcommand {\rrr}[1]{(\ref{#1})}

\def \be{\begin{equation}}
\def \ee{\end{equation}}
\def \bt{\begin{theorem}}
\def \et{\end{theorem}}
\def \bc{\begin{corollary}}
\def \ec{\end{corollary}}
\def \bea{\begin{eqnarray}}
\def \eea{\end{eqnarray}}
\def \bas{\begin{eqnarray*}}
\def \eas{\end{eqnarray*}}

\def \noi{\noindent}

\def \gg{\gamma}

\def \la{\lambda}

\def \pp{\phi}

\def \vski{\vspace{12pt}}

\def \VV{{\cal V}}

\def \({\left(}
\def \){\right)}

\def \nn{\nonumber}

\def \bc{\begin{center} }
\def \ec{\end{center} }
\def \bs{\begin{slide} }
\def \es{\end{slide} }

\def\square{{\vcenter{\vbox{\hrule height.3pt
         \hbox{\vrule width.3pt height5pt \kern5pt
            \vrule width.3pt}
         \hrule height.3pt}}}}
\def\qed{{\hfill $\square$ \bigskip}}

\newcounter{cccases}
\newcommand{\ccases}[1]{\begingroup \refstepcounter{cccases} {\bf \fontsize{14}{16}\selectfont Case \thecccases }  \label{#1}\endgroup}

\begin{document}

\title{On electric resistances for distance-regular graphs}

 \author{Jack H. Koolen$^{\,\rm 1}$,~Greg Markowsky$^{\,\rm 2}$ and Jongyook Park$^{\,\rm 1}$\\
{\small {\tt koolen@postech.ac.kr} ~~ {\tt gmarkowsky@gmail.com} ~~ {\tt jongyook@postech.ac.kr}}\\
{\footnotesize{$^{\rm 1}$Department of Mathematics,  POSTECH, Pohang 790-785, South Korea}}\\
{\footnotesize{$^{\rm 2}$Pohang Mathematics Institute,  POSTECH, Pohang 790-785, South Korea}}}
\maketitle

\begin{abstract}
We investigate the behavior of electric potentials on distance-regular graphs, and extend some results of a prior paper, \cite{markool}. Our main result, Theorem \ref{thispap} below, shows(together with Corollary \ref{last}) that if distance is measured by the electric resistance between points then all points are close to being equidistant on a distance-regular graph with large valency. A number of auxiliary results are also presented.
\end{abstract}

\section{Introduction}

In \cite{biggs2}(or see also \cite{biggs}), Biggs laid the foundation for studying electric resistance on distance-regular graphs. Such graphs were shown to be very natural objects in the study of potential theory, due to the ease of making explicit calculations. In the same paper, Biggs made the following conjecture, which has recently been proved in \cite{markool}.

\bt \lll{bigguy}
Let G be a distance-regular graph with degree larger than 2 and diameter D. If $d_j$ is the electric resistance between
any two vertices of distance $j$, then

\be
\max_j d_j = d_D \leq K d_1
\ee

\noi where $K=1+\frac{94}{101} \approx 1.931$. Equality
holds only in the case of the Biggs-Smith graph.
\et

\noi In this paper we prove several results which extend this theorem, the primary of which is the following.

\bt \lll{latcor}
Let G be a distance-regular graph with degree larger than 2 and diameter $D \geq 3$. Then

\be
\max_j d_j = d_D < (1+\frac{6}{k}) d_1
\ee

\noi where $k$ is the degree of $G$.
\et

This will be seen to be a corollary of Theorem \ref{thispap} below. Note that $d_1 = \min_{1 \leq j \leq D} d_j$. Theorem \ref{latcor} therefore shows, essentially, that under the electric resistance metric all points are more or less equidistant in distance-regular graphs of large degree(other than the strongly regular graphs, which are handled in Section \ref{srgraphs}). We will proceed as follows. Section \ref{drgprop} gives the definition of distance-regular graphs, while Section \ref{prevkno} gives known properties of electric resistance on this class of graphs. Section \ref{newres} states the results of this paper, while Section \ref{test} gives a proof of our main result. Finally, Section \ref{srgraphs} discusses the case in which the underlying graph is strongly regular and Section \ref{conc} gives a few concluding remarks.

\section{Distance-regular graphs} \lll{drgprop}

All the graphs considered in this paper are finite, undirected and
simple (for unexplained terminology and more details, see for example \cite{drgraphs}). Let $G$ be a
connected graph and let $\VV=\VV(G)$ be the vertex set of $G$. The distance $d(x,y)$ between
any two vertices $x,y$ of $G$
is the length of a shortest path between $x$ and $y$ in $G$. The diameter of $G$ is the maximal distance
occurring in $G$ and we will denote this by $D = D(G)$.  
For a vertex $x \in \VV(G)$, define $K_i(x)$ to be the set of
vertices which are at distance $i$ from $x~(0\le i\le
D)$. In addition, define $K_{-1}(x):=\emptyset$ and $K_{D+1}(x)
:= \emptyset$. We write $x\sim_{G} y$ or simply $x\sim y$ if two vertices $x$ and $y$ are adjacent in $G$. A connected graph $G$ with diameter $D$ is called
{\em distance-regular} if there are integers $b_i,c_i$ $(0 \le i
\leq D)$ such that for any two vertices $x,y \in \VV(G)$ with $d(x,y)=i$, there are precisely $c_i$
neighbors of $y$ in
$K_{i-1}(x)$ and $b_i$ neighbors of $y$ in $K_{i+1}(x)$
(cf. \cite[p.126]{drgraphs}). In particular, distance-regular graph $G$ is regular with valency
$k := b_0$ and we define $a_i:=k-b_i-c_i$ for notational convenience. 
The numbers $a_i$, $b_{i}$ and $c_i~(0\leq i\leq D)$ are called the {\em
intersection numbers} of $G$. Note that $b_D=c_0=a_0=0$, $b_0 = k$ and $c_1=1$.
The intersection numbers of a distance-regular graph $G$ with diameter $D$ and valency $k$ satisfy
(cf. \cite[Proposition 4.1.6]{drgraphs})\\

(i) $k=b_0> b_1\geq \cdots \geq b_{D-1}$;\\
(ii) $1=c_1\leq c_2\leq \cdots \leq c_{D}$;\\
(iii) $b_i\ge c_j$ \mbox{ if }$i+j\le D$.\\

Moreover, if we fix a vertex $x$ of $G$, then $|K_i(x)|$ does not depend on the
choice of $x$ as $c_{i+1} |K_{i+1}(x)| =
b_i |K_i(x)|$ holds for $i =1, 2, \ldots D-1$. We will therefore write $K_i$ instead of $K_i(x)$ where convenient. In the next section, it will be shown that the resistance between any two vertices
of $G$ can be calculated explicitly using only the intersection array, so that the proof can be conducted using only
the known properties of the array.

\section{Electric resistance on distance-regular graphs} \lll{prevkno}

Henceforth let $G$ be a distance-regular graph with $n$ vertices, degree $k \geq 3$, and diameter $D$. To calculate the resistance between any two vertices we use Ohm's Law, which states that

\be \label{}
V=IR
\ee

\noi where $V$ represents a difference in voltage(or potential), $I$ represents current, and $R$ represents resistance.
That is, we imagine that our graph is a circuit where each edge is a wire with resistance 1. We attach a battery of
voltage $V$ to
two distinct vertices $u$ and $v$, producing a current through the graph. The resistance between $u$ and $v$
is then $V$
divided by the current produced. The current flowing through the circuit can be determined by calculating the voltage at
each point on the graph, then summing the currents flowing from $u$, say, to all vertices adjacent to $u$. Calculating the voltage at each point is thereby seen to be an important problem. A function $f$ on the vertex set $\VV$ is harmonic at a point $z\in \VV$ if $f(z)$ is the average of neighboring values of $f$, that is

\be \label{}
\sum_{x\sim z} (f(x)-f(z)) = 0
\ee

The voltage
function on $\VV$ can be characterized as the unique function which is harmonic on $\VV-\{u,v\}$ having the prescribed values
on $u$ and $v$. For our purposes, on the distance-regular graph $G$, we will first suppose that $u$ and $v$ are adjacent.
It is easy to see that, for any vertex $z$, $|d(u,z)-d(v,z)| \leq 1$, where $d$ denotes the ordinary graph-theoretic
distance. Thus, any $z$ must be contained in a unique set of one of the following forms:

\bea \label{}
&& K_i^i  = \{x: d(u,x)=i \mbox{ and } d(v,x)=i \}
\\ \nn && K_i^{i+1} = \{x: d(u,x)=i+1 \mbox{ and } d(v,x)=i \}
\\ \nn && K_{i+1}^i = \{x: d(u,x)=i \mbox{ and } d(v,x)=i+1 \}
\eea

\noi Suppose that $(b_0,b_1,\ldots ,b_{D-1};c_1,c_2,\ldots,c_D)$ is the intersection array of $G$. For $0 \leq i \leq D-1$
define the numbers $\phi_i$ recursively by

\bea \label{smile}
&& \pp_0=n-1
\\ \nn && \pp_i = \frac{c_i\pp_{i-1}-k}{b_i}
\eea

It can be shown(see \cite{biggs2} or \cite{markool}) that $\phi_0,\phi_1, \ldots, \phi_{D-1}$ is a strictly decreasing sequence. The explicit value of $\phi_i$ is given by the following equation, first stated by Biggs in \cite{biggs2}.

\be \label{nut}
\phi_i = k\Big(\frac{1}{c_{i+1}} + \frac{b_{i+1}}{c_{i+1}c_{i+2}} + \ldots + \frac{b_{i+1} \ldots b_{D-1}}{c_{i+1} \ldots c_{D}} \Big)
\ee

We then have the following fundamental proposition from \cite{biggs2}.

\begin{proposition} \label{vp}
The function $f$ defined on $\VV$ by

\bea \label{}
&& f(u) = -f(v) = \pp_0
\\ \nn && f(z)= 0 \mbox{ for } x \in K_{i}^i
\\ \nn && f(z)= \pp_i \mbox{ for } x \in K_{i+1}^i
\\ \nn && f(z)= -\pp_i \mbox{ for } x \in K^{i+1}_i
\eea

is harmonic on $\VV-\{u,v\}$.
\end{proposition}

This leads to the following, also given in \cite{biggs2}.

\begin{corollary} \label{tool}
The resistance between two vertices of distance $j$ in a graph is given by
\be \label{}
d_j = \frac{2\sum_{0\leq i<j}\pp_i}{nk}
\ee
\end{corollary}

Theorem \ref{bigguy} and Corollary \ref{tool} imply that

\be
d_{D} \leq \frac{4\pp_0}{nk} < \frac{4}{k}
\ee

so that the maximal resistance between any two points on a distance-regular graph goes to 0 as the degree of the graph goes to infinity. Note also that Theorem \ref{bigguy} implies the following.

\be \label{}
\frac{d_D}{d_1} = \frac{\phi_0 + \ldots + \phi_{D-1}}{\phi_{0}} \leq \frac{195}{101}
\ee

with equality holding only in the case of Biggs-Smith graph.

\section{Results on $\phi_i$ for $i \geq 1$} \lll{newres}

The discussion in the previous section indicates that understanding the behavior of the $\phi_i$'s is crucial for understanding the theory of electric resistance on distance-regular graphs. In this section, we present some further results. Theorem \ref{bigguy} shows that $\phi_0$ bounds the sum of the terms which follow. Similar behavior is exhibited by the ensuing $\phi_i$'s, as the following theorem indicates.

\bt \lll{iii}
For any $m \geq 0$,

\be \label{}
\phi_{m+1} + \ldots + \phi_{D-1} < (3m+3)\phi_m
\ee
\et

{\bf Proof:} If $D \leq 4(m+1)$, this is trivial due to the strict monotonicity of the $\phi_i$'s. Suppose $D>4(m+1)$. By Lemma 5.2 of \cite{hbk}, $l \geq m+2$, where $l$ is the unique integer which satisfies $\frac{b_{l-1}}{c_{l-1}}\geq 2$ but $\frac{b_l}{c_l}<2$. Again applying Lemma 5.2 of \cite{hbk}, we see that $D \leq 4l$. Using the fact that

\be \lll{comm}
\frac{\phi_i}{\phi_{i-1}} < \frac{c_i}{b_i}
\ee

and the monotonicity of the $\phi_i$'s, we have

\bea \label{}
&&
\frac{\phi_{m+1} + \ldots + \phi_{D-1}}{\phi_m} < \frac{1}{2} + \frac{1}{4} + \ldots + \frac{1}{2^{l-m-1}} + \frac{1}{2^{l-m-1}}(3l)
\\ \nn && \hspace{1cm} < 1 + \frac{3l}{2^{l-m-1}}
\eea

The last expression is decreasing in $l$, so we can bound it by plugging in $l=m+2$ to get

\be \label{}
\frac{\phi_{m+1} + \ldots + \phi_{D-1}}{\phi_m} < 1 + \frac{3(m+2)}{2}<3m+3
\ee \qed

If we take $m=1$ we get $\phi_{2} + \ldots + \phi_{D-1} < 6 \phi_1$. However, this can be improved quite a bit as is shown by the following theorem, which is the main result of this paper.

\bt \lll{thispap}
\be
\phi_{2} + \ldots + \phi_{D-1} \leq\phi_1
\ee

with equality holding only in the case of the dodecahedron.
\et

The proof of this theorem relies on analyzing a number of cases, and is given in Section \ref{test}. This theorem yields an interesting consequence when combined with Corollary \ref{tool}. Let us consider the ratio of $d_1$ with the maximal possible resistance on $G$, $d_D$. This is given by

\be \label{}
\frac{d_D}{d_1} = \frac{\phi_0 + \ldots + \phi_{D-1}}{\phi_{0}}
\ee

Theorem \ref{bigguy} implies that this ratio is bounded by 2. Theorem \ref{thispap} improves this bound, for we see that

\be \label{}
\frac{d_D}{d_1} \leq \frac{\phi_0 + 2\phi_1}{\phi_{0}}
\ee

If $G$ is not strongly regular(diameter 2), then $b_1 \geq k/3$ by Lemma 3 of \cite{parkool2}, and since $c_1$ is necessarily 1 we see that $\phi_1 < \frac{3}{k} \phi_0$. This allows us to obtain the following corollary, referred to earlier as Theorem \ref{latcor}.

\begin{corollary} \label{hhh22}
\be
\frac{d_D}{d_1} \leq 1+\frac{6}{k}
\ee
\end{corollary}

This shows that for large $k$, all points become nearly equidistant when measured with respect to the resistance metric. We remark further that it has been shown that $b_1 \geq k/2$ in all but a small number of known families of distance-regular graphs; see Theorem 11 in \cite{parkool2}. We may therefore replace the $1+\frac{6}{k}$ by $1+\frac{4}{k}$ in many cases if desired. We would like to thank Professor Tommy Jensen for posing the question to us whether something like Corollary \ref{hhh22} might be true.

\ind In light of Theorems \ref{bigguy} and \ref{thispap} it may be tempting to hope that $\phi_{m+1} + \ldots + \phi_{D-1} \leq \phi_m$ for all $m$. However, the Biggs-Smith graph, with intersection array $(3,2,2,2,1,1,1;1,1,1,1,1,1,3)$, yields the following potentials.

\be \label{}
\phi_0=101, \phi_1=49, \phi_2=23, \phi_3=10, \phi_4=7, \phi_5=4, \phi_6=1
\ee

Note that $\phi_4 + \phi_5 + \phi_6 > \phi_3$. Nonetheless, we do conjecture the following strengthening of Theorem \ref{iii}.

\begin{conjecture}

There is a universal constant $K$ such that

\be
\phi_{m+1} + \ldots + \phi_{D-1} \leq K \phi_m
\ee

for all $m$ and all distance-regular graphs.
\end{conjecture}

\vski

It is an educated guess that $2$ might work for $K$, but we do not at this stage have a proof.

\section{Proof of Theorem \ref{thispap}} \lll{test}

As in the proof of Theorem \ref{bigguy} in \cite{markool}, we will consider a number of special cases.

\vski

\ccases{smallD}: $D \leq 3$.

There is nothing to prove when $D=2$, and the case $D=3$ is trivial, since $\phi_2<\phi_1$. \hfill $\triangle$

\vski

\ccases{deg3}: $k=3$.

It is known(see \cite{drgraphs}, Theorem 7.5.1) that the only distance-regular graphs of degree 3 with
diameter greater than 3 are given by the intersection arrays below, and which give rise to the resistances below:

\begin{tabular}{ l c c r }
Name & Vertices & Intersection array & $\frac{\phi_2+ \ldots \phi_{D-1}}{\phi_1}$ \\
\hline
Pappus graph & 18 & (3,2,2,1;1,1,2,3)  &      0.428571 \\
Coxeter graph & 28 & (3,2,2,1;1,1,1,2)   &     0.5 \\
Tutte's 8-cage & 30 & (3,2,2,2;1,1,1,3)    &    0.461538 \\
Dodecahedron & 20 & (3,2,1,1,1;1,1,1,2,3) &   1.0 \\
Desargues graph & 20 & (3,2,2,1,1;1,1,2,2,3)  &  0.6875 \\
Tutte's 12-cage & 126 & (3,2,2,2,2,2;1,1,1,1,1,3)  &     0.786885 \\
Biggs-Smith graph & 102 & (3,2,2,2,1,1,1;1,1,1,1,1,1,3) &  0.918367 \\
Foster graph & 90 & (3,2,2,2,2,1,1,1;1,1,1,1,2,2,2,3)  &      0.854651 \\
\end{tabular}

\vski

\ccases{deg4}: $k=4$.

It is known(see \cite{broukool}) that the only distance-regular graphs of degree 4 with diameter greater than 3
are given by the intersection arrays below, and which give rise to the resistances below:

\noi \begin{tabular}{ l c c r }
Name & Vertices & Intersection array & $\frac{\phi_2+ \ldots \phi_{D-1}}{\phi_1}$ \\
\hline
4-cube &  16 & (4,3,2,1;1,2,3,4) &       0.727273\\
Flag graph of $PG(2,2)$ & 21 & (4,2,2;1,1,2) &   0.25\\
Incidence graph of $PG(2,3)$ & 26 & (4,3,3;1,1,4)  &  0.142857\\
Incidence graph of $AG(2,4)$-p.c. & 32 & (4,3,3,1;1,1,3,4)  &      0.296296\\
Odd graph $O_4$ & 35 & (4,3,3;1,1,2)  &  0.2\\
Flag graph of $GQ(2,2)$ & 45 & (4,2,2,2;1,1,1,2)  &      0.5\\
Doubled odd graph & 70& (4,3,3,2,2,1,1;1,1,2,2,3,3,4)&    0.661538\\
Incidence graph of $GQ(3,3)$ & 80 & (4,3,3,3;1,1,1,4) &       0.32\\
Flag graph of $GH(2,2)$ & 189 & (4,2,2,2,2,2;1,1,1,1,1,2) &      0.804348\\
Incidence graph of $GH(3,3)$ & 728 & (4,3,3,3,3,3;1,1,1,1,1,4) &      0.46473\\
\end{tabular}

\vski

\ccases{niceb2c2}: $\frac{b_2}{c_2} \geq 10$.

\vski

Applying Theorem \ref{iii}, we have $\phi_2 + \ldots + \phi_{D-1} \leq 10 \phi_2$, so the result follows on noting $\phi_2 < \frac{c_2}{b_2} \phi_1$. \hfill $\triangle$

\vski

\ccases{smallD4}: $D = 4, c_3>1$.

It is a consequence of \cite[Theorem 5.4.1]{drgraphs} that $c_3 \geq \frac{3}{2}c_2$. Using \rrr{nut}, we have

\be \lll{fff}
\frac{c_2c_3c_4}{k}(\phi_1-\phi_2-\phi_3) = c_3c_4 + b_2c_4 + b_2b_3 -c_2c_4 - c_2b_3 - c_2c_3
\ee

But $b_2 \geq c_2$ by Property $(iii)$ of Section \ref{drgprop}, and thus $b_2c_4 \geq c_2c_3$, $b_2b_3 \geq c_2b_3$. Since $c_3c_4 \geq \frac{3}{2}c_2c_4$, \rrr{fff} is strictly positive. \hfill $\triangle$

\vski

\ccases{smallD5}: $D = 5, c_3>1$.

\vski

The following lemma will be useful in this and subsequent cases.

\begin{lemma} \lll{mj}
Suppose $c_3\geq2$ and $D\geq5$. Then $c_3 \geq \frac{3}{2}c_2$, and either $c_4\geq2c_2$ or $c_5\geq c_2+c_3$. In both cases we have $c_5 \geq 2 c_2$, and if $c_5\geq c_2+c_3$ then $c_5\geq \frac{5}{2}c_2$.
\end{lemma}

{\bf Proof:} If $c_2\geq2$, then $c_3\geq\frac{3}{2}c_2$ holds by \cite[Theorem 5.4.1]{drgraphs}, while if $c_2=1$ then $c_3 \geq 2 = 2c_2$, proving that $c_3 \geq \frac{3}{2}c_2$ in all cases. Let $j\geq4$ be the minimal integer such that $c_{j-1}<c_j$. Then $j=4$ or $5$ by \cite[Theorem 1.1]{hbk}. If $j=4$(i.e., $c_3<c_4$), then $c_4\geq c_2+c_{4-2}=2c_2$ holds by \cite[Proposition 1 (ii)]{koolen}. Similarly, $c_5\geq c_2+c_{5-2}=c_2+c_3$ holds when $j=5$. This proves the first statement, and the last statement is a simple consequence of the first. This proves the lemma. \hfill $\lozenge$

\vski

Now, if $\frac{b_2}{c_2} < \frac{3}{2}$, Lemma \ref{mj} shows that $c_3>b_2$, contradicting Property $(iii)$ of Section \ref{drgprop}. Thus, we may assume $\frac{b_2}{c_2} \geq \frac{3}{2}$. Thus, $\phi_1-\phi_2 = \Big(\frac{b_2\phi_2+k}{c_2}\Big) - \phi_2 \geq \frac{k}{c_2} + \frac{1}{2}\phi_2$. We therefore have

\bea \label{}
&& \frac{c_3c_4c_5}{k}(\phi_1-\phi_2-\phi_3-\phi_{4})
\\ \nn && \hspace{1cm} \geq \frac{c_3c_4c_5}{k}(\frac{k}{c_2}+\frac{1}{2}\phi_2-\phi_3-\phi_{4})
\\ \nn && \hspace{1cm} \geq c_3c_4c_5(\frac{1}{c_2}+\frac{1}{2}(\frac{1}{c_3}+\frac{b_3}{c_3c_4} + \frac{b_3b_4}{c_3c_4c_5})-(\frac{1}{c_4}+\frac{b_4}{c_4c_5})-\frac{1}{c_5})
\\ \nn && \hspace{1cm} = (\frac{c_4}{c_2}-1)c_3c_5+\frac{1}{2}(c_4c_5+b_3c_5 + b_3b_4)-c_3c_4-c_3b_4
\eea

If $c_4\geq 2 c_2$, then $(\frac{c_4}{c_2}-1) \geq 1$, and we are reduced to showing

\be \lll{}
c_3c_5+\frac{1}{2}(c_4c_5+b_3c_5 + b_3b_4)>c_3c_4+c_3b_4
\ee

This is clear, since $c_3c_5 \geq c_3c_4$ and $b_3c_5 \geq c_3b_4$, while $c_4c_5 \geq c_3b_4$ if $b_4 \leq c_4$ whereas $b_3b_4 > c_3b_4$ if $b_4 > c_4$. Now, if $c_4<2c_2$ we still have $c_4 \geq \frac{3}{2}c_2$ by \cite[Theorem 5.4.1]{drgraphs}, and Lemma \ref{mj} and its proof show that $c_3=c_4$. We must therefore show

\be \lll{rent}
\frac{1}{2}(c_3c_5+c_4c_5+b_3c_5 + b_3b_4)>c_3c_4+c_3b_4 = c_3^2 + c_3b_4
\ee

Also by Lemma \ref{mj} we have $c_5\geq c_2+c_3$. We see that
\bea \lll{rent2} \nn
\\ \nn && \frac{1}{2}(c_3c_5+c_4c_5+b_3c_5 + b_3b_4)> \frac{1}{2}(2c_3^2+2c_2c_3+b_3c_5 + b_3b_4)
\\ \\ \nn && \hspace{1cm} > \frac{3}{2} c_3^2 + \frac{1}{2}(b_3c_5 + b_3b_4)
\eea

In all cases $b_3c_5 \geq c_3b_4$, and if $c_3 \geq b_3$ then $c_3^2 \geq c_3b_4$, whereas if $c_3< b_3$ then $b_3b_4 > c_3 b_4$. We conclude either way that

\be \label{rent3}
\frac{3}{2} c_3^2 + \frac{1}{2}(b_3c_5 + b_3b_4) \geq c_3^2 + c_3b_4
\ee

\rrr{rent2} and \rrr{rent3} together verify \rrr{rent}. \hfill $\triangle$

\vski

\ccases{niceD}: $6 \leq D \leq \frac{2b_2}{c_2}+2$, $c_3 > 1$.

\vski

Let $\aaa=D-4$, so that $\frac{b_2}{c_2} \geq 1+\frac{\aaa}{2}$, and

\be
\phi_1-\phi_2 = \Big(\frac{b_2\phi_2+k}{c_2}\Big) - \phi_2 \geq \frac{k}{c_2} + \frac{\aaa}{2}\phi_2
\ee

We will show $\frac{\phi_1-\phi_2-\cdots-\phi_{D-1}}{k}>0$. We have

\bea \label{}
&& \frac{\phi_1-\phi_2-\cdots-\phi_{D-1}}{k}
\\ \nn && \hspace{1cm} \geq \frac{\frac{k}{c_2}+\frac{\aaa}{2}\phi_2-\phi_3-\cdots-\phi_{D-1}}{k}
\\ \nn && \hspace{1cm} = \Bigg\{\frac{1}{c_2}+\frac{\alpha}{2}\times\frac{1}{c_3}-\frac{1}{c_4}-\frac{1}{c_5}-\cdots-\frac{1}{c_{\alpha+4}} \Bigg\}
\\ \nn && \hspace{1.5cm}+ \Bigg\{\frac{\alpha}{2}\Big(\frac{b_3}{c_3c_4}+\frac{b_3b_4}{c_3c_4c_5}+\cdots+\frac{b_3b_4\cdots b_{\alpha+3}}{c_3c_4c_5\cdots c_{\alpha+4}}\Big)
\\ \nn && \hspace{1.5cm} -\Big(\frac{b_4}{c_4c_5}+\frac{b_4b_5}{c_4c_5c_6}+\cdots+\frac{b_4b_5\cdots b_{\alpha+3}}{c_4c_5\cdots c_{\alpha+4}}\Big)
\\ \nn && \hspace{5cm} \vdots
\\ \nn && \hspace{1.5cm} - \Big(\frac{b_{\alpha+3}}{c_{\alpha+3}c_{\alpha+4}}\Big) \Bigg\}
\\ \nn && \hspace{1cm} =: \Big\{ I \Big\} + \Big\{ II \Big\}
\eea

We may now invoke the ensuing two lemmas. Lemma \ref{jimmy} below shows that $\{II\}$ is positive in the preceding expression, and Lemma \ref{river} shows that $\{I\}$ is nonnegative.

\begin{lemma} \lll{jimmy} For any $\gg$ with $2 \leq \gg \leq D-4$, the expression defined by $\{II\}$ is positive.
\end{lemma}

We will show this lemma by induction on $\gg$. If $\gg=2$, then $D \geq 6$, and thus $b_3\geq c_3$. We therefore know  $\frac{b_3b_4}{c_3c_4c_5}\geq\frac{b_4}{c_4c_5}$ and $\frac{b_3b_4b_5}{c_3c_4c_5c_6}\geq\frac{b_4b_5}{c_4c_5c_6}$. As $c_3\geq2$, we have $|\{i|c_i=c_3\}|\leq2$ (\cite[Theorem 1.1]{hbk}), which implies $c_5>c_3$. We see that $\frac{b_3}{c_3c_4}>\frac{b_5}{c_5c_6}$, as $b_3\geq b_5$ and $c_6\geq c_4$. This proves the lemma for $\gg=2$.

\vski

Suppose now that the lemma is true for $\gg$. Using the induction hypothesis, in order to obtain the result for $\gg+1$ we need to show

\bea \label{cov}
&& \frac{1}{2}(\frac{b_3}{c_3c_4}+\frac{b_3b_4}{c_3c_4c_5}+\cdots+\frac{b_3b_4\cdots b_{\gg+3}}{c_3c_4\cdots c_{\gg+3}c_{\gg+4}})+\frac{\gg+1}{2}\times\frac{b_3b_4\cdots b_{\gg+4}}{c_3c_4\cdots c_{\gg+4}c_{\gg+5}}
\\ \nn && \hspace{1cm} \geq\frac{b_4b_5\cdots b_{\gg+3}b_{\gg+4}}{c_4c_5\cdots c_{\gg+3}c_{\gg+4}c_{\gg+5}}+\frac{b_5b_6\cdots b_{\gg+3}b_{\gg+4}}{c_5c_6\cdots c_{\gg+3}c_{\gg+4}c_{\gg+5}}+\cdots+\frac{b_{\gg+4}}{c_{\gg+4}c_{\gg+5}}
\eea

Note that

\be \lll{spy}
\frac{b_ib_{i+1}\cdots b_j}{c_ic_{i+1}\cdots c_j}=\frac{b_i}{c_j}\times\frac{b_{i+1}}{c_{j-1}}\times\cdots\times\frac{b_j}{c_i}\geq 1
\ee

if $i+j\leq D$. This will be an important fact for us. We first suppose that $\gg+1$ is an even integer. We can check easily that

\bea
&& \nn \frac{\gg+1}{2}\times\frac{b_3b_4\cdots b_{\gg+4}}{c_3c_4\cdots c_{\gg+4}c_{\gg+5}} \geq\frac{b_4b_5\cdots b_{\gg+3}b_{\gg+4}}{c_4c_5\cdots c_{\gg+3}c_{\gg+4}c_{\gg+5}}+\frac{b_5b_6\cdots b_{\gg+3}b_{\gg+4}}{c_5c_6\cdots c_{\gg+3}c_{\gg+4}c_{\gg+5}} +
\\ \nn && \hspace{1cm} \cdots+\frac{b_{(\gg+1)/2+3}b_{(\gg+1)/2+4}\cdots b_{\gg+3}b_{\gg+4}}{c_{(\gg+1)/2+3}c_{(\gg+1)/2+4}\cdots c_{\gg+3}c_{\gg+4}c_{\gg+5}}
\eea

due to \rrr{spy} and the fact that $3+(\frac{\gg+1}{2}+2) \leq D$. Also, we may show that

\bea
&& \nn \frac{1}{2}(\frac{b_3b_4\cdots b_{\gg+3-2i}}{c_3c_4\cdots c_{\gg+4-2i}}+\frac{b_3b_4\cdots b_{\gg+2-2i}}{c_3c_4\cdots c_{\gg+3-2i}})
\\ \nn && \hspace{.2cm} \geq \frac{1}{2}(\frac{b_3\cdots b_{(\gg+1)/2+2-i}b_{(\gg+1)/2+4+i}\cdots b_{n+3}}{c_3\cdots c_{(\gg+1)/2+2-i}c_{(\gg+1)/2+4+i}\cdots c_{\gg+4}}+\frac{b_3\cdots b_{(\gg+1)/2+1-i}b_{(\gg+1)/2+4+i}\cdots b_{n+3}}{c_3\cdots c_{(\gg+1)/2+1-i}c_{(\gg+1)/2+4+i}\cdots c_{\gg+4}})
\\ \nn && \hspace{.2cm} \geq\frac{b_{(\gg+1)/2+4+i}\cdots b_{\gg+3}b_{\gg+4}}{c_{(\gg+1)/2+4+i}\cdots c_{\gg+3}c_{\gg+4}c_{\gg+5}}
\eea

holds for $i=0,1,\ldots,\frac{n-1}{2}$. To see this, note that the first inequality is obtained by merely shifting some of the indices, resulting in smaller quotients, and the second inequality is due again to \rrr{spy} and the fact that $3+(\gg+1)/2+2-i \leq D$. Thus, \rrr{cov} is true when $\gg+1$ is an even integer.

\vski

Now let us suppose that $\gg+1$ is an odd integer. By the same arguments as above, we have

\bea
&& \frac{\gg}{2}\times\frac{b_3b_4\cdots b_{\gg+4}}{c_3c_4\cdots c_{\gg+4}c_{\gg+5}}\geq\frac{b_4b_5\cdots b_{\gg+3}b_{\gg+4}}{c_4c_5\cdots c_{\gg+3}c_{\gg+4}c_{\gg+5}}+\frac{b_5b_6\cdots b_{\gg+3}b_{\gg+4}}{c_5c_6\cdots c_{\gg+3}c_{\gg+4}c_{\gg+5}}+
\\ \nn && \hspace{1cm} \cdots+\frac{b_{\gg/2+3}b_{\gg/2+4}\cdots b_{\gg+3}b_{\gg+4}}{c_{\gg/2+3}c_{\gg/2+4}\cdots c_{\gg+3}c_{\gg+4}c_{\gg+5}}
\eea

and

\be
\frac{1}{2}(\frac{b_3b_4\cdots b_{\gg+4-2i}}{c_3c_4\cdots c_{\gg+5-2i}}+\frac{b_3b_4\cdots b_{\gg+3-2i}}{c_3c_4\cdots c_{\gg+4-2i}})\geq\frac{b_{\gg/2+4+i}\cdots b_{\gg+3}b_{\gg+4}}{c_{\gg/2+4+i}\cdots c_{\gg+3}c_{\gg+4}c_{\gg+5}}
\ee

for $i=0,1,\ldots,\frac{n}{2}$. Thus, \rrr{cov} is true if $n+1$ is an odd integer. This completes the induction and the proof of the lemma. \hfill $\lozenge$

\begin{lemma} \lll{river} For any $\gg$ with $0 \leq \gg \leq D-4$, we have
\be \lll{time}
\frac{1}{c_2}+\frac{\gg}{2}\times\frac{1}{c_3}\geq\frac{1}{c_4}+\frac{1}{c_5}+\cdots+\frac{1}{c_{\alpha+4}}
\ee
\end{lemma}

To prove this, let us consider several cases. We begin with $\gg=1$, as the case $\gg=0$ is trivial. Note that we are still assuming $c_3>1$.

\vski

$(i)$ For $\gg=1$, by Lemma \ref{mj} we have $c_5 \geq 2c_2$. As a consequence of Theorem 5.4.1 of \cite{drgraphs} and the fact that $c_3>1$, $c_2 \leq \frac{2}{3}c_3$, whence $\frac{1}{3} \times \frac{1}{c_2} \geq \frac{1}{2}\times \frac{1}{c_3}$. Thus,

\be
\frac{1}{c_2} + \frac{1}{2}\times\frac{1}{c_3} \geq \frac{2}{3} \times \frac{1}{c_2} + \frac{1}{c_3} \geq \frac{2}{3} \times 2 \times \frac{1}{c_5} + \frac{1}{c_3} > \frac{1}{c_4} + \frac{1}{c_5}.
\ee

\vski

$(ii)$ For $\gg=2$, we know that $c_5\geq2c_2$ holds by Lemma \ref{mj}. Thus, $\frac{1}{c_2} \geq\frac{2}{c_5}\geq\frac{1}{c_5}+\frac{1}{c_6}$, and the result follows by noting $c_3\leq c_4$.

\vski

$(iii)$ For $\gg=3$, we first assume $c_5>c_4$, which implies $c_5\geq c_2+c_3$ by \cite[Proposition 1 (ii)]{koolen}. If $c_3\geq2c_2$, then $c_5\geq3c_2$ holds, hence $\frac{1}{c_2}\geq\frac{3}{c_5}\geq\frac{1}{c_5}+\frac{1}{c_6}+\frac{1}{c_7}$. The inequality therefore holds, as $\frac{3}{2}\times\frac{1}{c_3}\geq\frac{1}{c_4}$. So we may assume $c_3<2c_2$. i.e., $\frac{1}{c_3}>\frac{1}{2}\times\frac{1}{c_2}$. Note that $c_5\geq\frac{5}{2}c_2$ by Lemma \ref{mj}. This implies $\frac{1}{c_2}+\frac{1}{2}\times\frac{1}{c_3}>\frac{5}{4}\times\frac{1}{c_2}>\frac{6}{5}\times\frac{1}{c_2}\geq\frac{1}{c_5}+\frac{1}{c_6}+\frac{1}{c_7}$. The inequality therefore holds, as $c_3\leq c_4$. Now suppose $c_5=c_4$. The proof of Lemma \ref{mj} implies that $c_4 \geq 2c_2$. By \cite[Theorem 1.1]{hbk}, we know that $|\{i|c_i=c_4\}|\leq 3$, which implies that either $c_6>c_5$ or $c_7>c_6$. Both cases imply $c_7\geq2c_3$ by \cite[Proposition 1 (ii)]{koolen}, and thus $\frac{1}{2}\times\frac{1}{c_3}\geq\frac{1}{c_7}$. We therefore have $\frac{1}{c_2} \geq \frac{1}{c_4}+\frac{1}{c_5}$ and $\frac{3}{2}\times\frac{1}{c_3} \geq \frac{1}{c_4} + \frac{1}{c_7}$. The result follows.

\vski

$(iv)$ For $\gg=4$, suppose $c_i>c_{i-1}$ holds for some $i\in\{6,7,8\}$. Then, by \cite[Proposition 1 (ii)]{koolen}, $c_i\geq2c_3$. This is enough to give the result, for in all cases $2c_3\leq c_8$, and thus we may replace $\frac{1}{c_8}$ by $\frac{1}{2} \times \frac{1}{c_3}$ in the right side of \rrr{time} and appeal to $(iii)$. We may therefore assume $c_5=c_6=c_7=c_8$. By \cite[Theorem 1.1]{hbk}, we know that $|\{i|c_i=c_4\}|\leq 3$, and thus $c_5>c_4$. As $\frac{1}{c_3}\geq\frac{1}{c_4}$, we want to show $\frac{1}{c_2}+\frac{1}{c_3}\geq\frac{1}{c_5}+\frac{1}{c_6}+\frac{1}{c_7}+\frac{1}{c_8}=\frac{4}{c_5}$. To show this, let us assume $\frac{1}{c_2}+\frac{1}{c_3}<\frac{4}{c_5}$ and derive a contradiction. Let $c_5=\beta c_2$ (We know $\beta \geq 2$ by Lemma \ref{mj}). Then $c_5\geq c_2+c_3$ by \cite[Proposition 1 (ii)]{koolen}, and thus $c_3\leq(\beta-1)c_2$. We see that $\frac{1}{\beta-1}\times\frac{1}{c_2}\leq\frac{1}{c_3}$ holds. As $\frac{1}{c_3}<\frac{4}{c_5}-\frac{1}{c_2}=(\frac{4}{\beta}-1)\times\frac{1}{c_2}$, we have $\frac{1}{\beta-1}<\frac{4}{\beta}-1$, which implies $(\beta-2)^2<0$. Clearly, this is a contradiction, and therefore the desired inequality holds.

\vski

$(v)$ For $\gg\geq5$, we know $c_9\geq2c_3$ by \cite[Proposition 1 (ii)]{koolen} and \cite[Theorem 1.1]{hbk}. Thus, $\frac{1}{2}\times\frac{1}{c_3}\geq\frac{1}{c_i}$ for all integers $i\geq9$. The inequality therefore holds by induction, using $(iv)$ as our initial case. This completes the proof of the lemma, and Case \ref{niceD} is now complete. \hfill $\triangle$

\vski

\ccases{c2>1}: $c_2 > 1$, $D \geq 6$, $G$ contains a quadrangle.

\vski

The following  is Lemma 11 in \cite{parkool}.

\begin{lemma} \label{}
Assume $D\geq4$ and $c_2\geq2$. If $G$ contains a quadrangle and $\frac{b_2}{c_2}<\frac{\alpha}{2}$ holds for some integer $\alpha\geq2$, then $D \leq \alpha+1$.
\end{lemma}

\vski

It was shown in Case \ref{smallD5} that since $D \geq 5$, $\aaa \geq 3$. The lemma implies that $D \leq \frac{2b_2}{c_2}+2$, and we may therefore appeal to Case \ref{niceD}. \hfill $\triangle$

\vski

\ccases{c2>1}: $c_2 > 1$, $D \geq 6$, $G$ does not contain a quadrangle.

\vski

$G$ is what is called a {\it Terwilliger graph}. By \cite[Corollary 1.16.6]{drgraphs} there is no distance-regular Terwilliger graph with $D \geq 6$ and $k<50(c_2-1)$, so we may
assume $50(c_2-1)\leq k$. Also, we may assume $\frac{b_2}{c_2}<10$ due to Case \ref{niceb2c2}. If $a_2+c_2\geq\frac{1}{2}k$, then the second largest eigenvalue $\theta_1$ of $G$ is bigger than $\frac{b_1}{2}-1$, as $\theta_1>a_2+c_2\geq\frac{1}{2}k>\frac{b_1}{2}-1$. As $\theta_1>\frac{b_1}{2}-1$, there is no distance-regular Terwilliger graph satisfying $D \geq 6$ by Proposition 9 in \cite{parkool}. Thus, we may assume $k-b_2=a_2+c_2<\frac{1}{2}k$. i.e., $b_2>\frac{1}{2}k$. Then the inequality $20c_2> 2b_2>k\geq50c_2-50$ holds, and this inequality gives $c_2<\frac{5}{3}$, which contradicts $c_2\geq2$. \hfill $\triangle$

\vski

\ccases{c2c3=1}: $c_2 = c_3 = 1$, $b_2 \geq 4$.

\vski

The following lemma gives a number of key calculations required for this and ensuing cases. The quantities $j=\inf\{i:c_i \geq b_i\}$ and $h=\inf \{i:c_i > b_i\}$ will be fundamental in all that follows. Clearly $h \geq j$. From this point on we will define $\phi_i = 0$ for $i \geq D$; this will allow us to state a number of facts more simply.

\begin{lemma} \lll{bad}

The following relations hold.

\vski

i) $D-h \leq j-1$

\vski

ii) $0 \leq h-j \leq j-1$, unless $c_j=1$, in which case we have $0 \leq h-j \leq j$.

\vski

iii) $D \leq 3j-2$, unless $c_j=1$, in which case $D \leq 3j-1$.

\vski

iv) For $0 \leq i \leq h-j$, we have $\phi_{j-1} \geq \phi_{D-i}+\phi_{j-1+i} \geq \phi_{h+j-1-i} +\phi_{j-1+i}$.

\vski

v) $\phi_{j-1} \geq \phi_{2j-2} + \phi_{2j-1} > 2 \phi_{2j-1}$, unless $c_j=1$, in which case we still have $\phi_{j-1} \geq 2 \phi_{2j-1}$.

\vski

vi) $\phi_j + \ldots + \phi_{D-1} \leq (j-1) \phi_{j-1}$ unless $c_j=1$, in which case we still have $\phi_j + \ldots + \phi_{D-1} \leq (j-1/2) \phi_{j-1}$.

\vski

\end{lemma}

To prove $(i)$, note that we must have either $c_h > b_j$ or $c_j>b_h$, and in either case we must have $j+h>D$, whence $D+1 \leq j+h$ by Property $(iii)$ in Section \ref{drgprop}. For $(ii)$, let us note that Corollary 5.9.6 of \cite{drgraphs} immediately implies that $h-j\leq j$. Now assume $c_j>1$, and let $k=\inf\{i:c_i=c_j\}$. Theorem 1.1 of \cite{hbk} implies that $c_{2k-1}>c_k=c_j$, so $h \leq 2k-1$. Clearly $k \leq j$, and $(ii)$ follows. $(iii)$ follows easily by combining $(i)$ and $(ii)$. The second inequality in $(iv)$ is a simple consequence of $(i)$ and the monotonicity of the $\phi$'s, so we need only prove the first. The first inequality is trivial if $i=0$(since $\phi_D=0$), and for $i>0$ note that $\frac{b_j \ldots b_{j-1+i}}{c_j \ldots c_{j-1+i}} = 1$, so using \rrr{nut} we have

\bea \label{}
&& \frac{\phi_{j-1} - \phi_{j-1+i}}{k} = (\frac{1}{c_j} + \ldots +\frac{b_j\ldots b_{D-1}}{c_j\ldots c_{D}}) - (\frac{1}{c_{j+i}} + \ldots +\frac{b_{j+i}\ldots b_{D-1}}{c_{j+i}\ldots c_{D}})
\\ \nn && \hspace{1cm} = (\frac{1}{c_j} + \ldots +\frac{b_j\ldots b_{j+i-2}}{c_j\ldots c_{j+i-1}})
\\ \nn && \hspace{1cm} \geq (\frac{1}{c_{D-i+1}} + \ldots +\frac{b_{D-i+1}\ldots b_{D-1}}{c_{D-i+1}\ldots c_{D}}) = \frac{\phi_{D-i}}{k}
\eea

This establishes $(iv)$. Choosing $i=h-j$ in $(iv)$, together with $(ii)$ and $(i)$ imply $(v)$, since by $(ii)$ we have $h-1 \leq 2j-1$($h-1 \leq 2j-2$ when $c_j>1$) and by $(i)$ we have $D-h+j \leq 2j-1$ as well. To achieve $(vi)$, write

\bea \lll{}
&& \nn \phi_j + \ldots + \phi_{D-1} = (\phi_j + \ldots + \phi_{h-1}) + (\phi_h + \ldots + \phi_{h+j-2})
\\ \nn && \hspace{1cm} =: (I) + (II)
\eea

Suppose first that $c_j > 1$, so that $h-j \leq j-1$. There are $h-j$ terms in $(I)$, each of which by $(iv)$ can be paired with one of the $j-1$ terms in $(II)$ to form a total bounded by $\phi_{j-1}$. We arrive at $\phi_j + \ldots + \phi_{D-1} \leq (j-1)\phi_{j-1}$. If $c_j = 1$, then it can happen that $h-j=j$. In that case $(I)$ has $j$ terms, but each of the terms in $(II)$ can be paired by $(iv)$ to a term in $(I)$ to form a total bounded by $\phi_{j-1}$, leaving only the term $\phi_{h-1}=\phi_{2j-1}$ left. We therefore obtain the weaker bound $\phi_j + \ldots + \phi_{D-1} \leq (j-1)\phi_{j-1}+ \phi_{2j-1}$. Applying $(v)$ we obtain $\phi_j + \ldots + \phi_{D-1} \leq (j-1/2)\phi_{j-1}$, and $(vi)$ is proved. This completes the proof of the lemma. \hfill $\lozenge$

\vski

We may now dispose of Case \ref{c2c3=1}. Suppose first that $b_3=1$, so $j=3$. Using Lemma \ref{bad}$(vi)$ and \rrr{comm} we have

\be \lll{}
\frac{\phi_2+ \ldots +\phi_{D-1}}{\phi_1} < \frac{1}{b_2} + \frac{5}{2}\times\frac{1}{b_2}
\ee

As $b_2 \geq 4$, we see that this is bounded above by $\frac{7}{8}$. Now suppose $b_3 \geq 2$, so $j \geq 4$. If $j=4$, then again using Lemma \ref{bad}$(vi)$ and \rrr{comm} we have

\be \lll{}
\frac{\phi_2+ \ldots +\phi_{D-1}}{\phi_1} < \frac{1}{b_2} + \frac{1}{b_2b_3} + \frac{7}{2}\times\frac{1}{b_2b_3}
\ee

As $b_2\geq 4, b_3\geq 2$, this is bounded above by $\frac{13}{16}$. Now suppose $j \geq 5$. Let $\la = \frac{b_3-1}{b_3}$. For $4 \leq i < j$ we see $\frac{c_i}{b_i}\leq \la$, so once again using Lemma \ref{bad}$(vi)$ and \rrr{comm} we have

\bea \lll{jw}
&& \frac{\phi_2+ \ldots + \phi_{D-1}}{\phi_1} < \frac{1}{b_2} + \frac{1}{b_2b_3} + \frac{\la}{b_2b_3}+ \ldots + \frac{\la^{j-4}}{b_2b_3} + (j - 1/2)\frac{\la^{j-4}}{b_2b_3}
\\ \nn && \hspace{1cm} < \frac{2}{b_2} + (j - 1/2)\frac{\la^{j-4}}{b_2b_3}
\eea

The last inequality is due to the fact that $1+\la + \la^2 + \ldots = b_3$. As $b_2 \geq 4$, the result will follow if we show $f(j) := (j - 1/2)\frac{\la^{j-4}}{b_2b_3} \leq \frac{1}{b_2}$ for all $j \geq 5$. Fix $b_3$. Suppose that $j > b_3$. Then

\be \lll{}
\frac{f(j+1)}{f(j)} = \frac{b_3-1}{b_3} \times \frac{j+1/2}{j-1/2} \leq \frac{j-1}{j} \times \frac{j+1/2}{j-1/2} < 1
\ee

Thus, $f(j) > f(j+1)$, and this tells us that we may assume $j \leq b_3+1$. But in that case, \rrr{jw} is bounded by $\frac{1}{b_2}(2+\frac{b_3+1/2}{b_3} \times \la^{j-4}) < \frac{3}{b_2}<1$. This completes the proof of this case. \hfill $\triangle$

\vski

\ccases{c2=1}: $b_2 \geq 4, c_2 =1, c_3>1$.

\vski

We will use the $j$ and $h$ as defined in Case \ref{c2c3=1}. It is clear that we may use the better bounds of Lemma \ref{bad} in this case. Assume first that $D>16$. By Lemma 4.2 of \cite{hbk}, $b_4 \geq 2 c_5$, so that $\frac{b_4}{c_4} \geq 2$. Note that Lemma \ref{bad}$(iii)$ implies $j \geq 7$. Similarly to Case \ref{c2c3=1} we let $\la = \frac{b_2-1}{b_2}$, so that $\frac{c_i}{b_i} \leq \la$ for $5 \leq i < j$. We may then write

\bea \lll{jw2}
&& \frac{\phi_2+ \ldots + \phi_{D-1}}{\phi_1}
\\ \nn && \hspace{1cm} < \frac{1}{b_2} + \frac{1}{2b_2} + \frac{1}{4b_2} + \frac{\la}{4b_2}+ \frac{\la^2}{4b_2} +\ldots + \frac{\la^{j-5}}{4b_2} + (j - 1)\frac{\la^{j-5}}{4b_2}
\\ \nn && \hspace{1cm} < \frac{3}{2b_2} + \frac{1}{4} + (j - 1)\frac{\la^{j-5}}{4b_2}
\eea

where we have used Lemma \ref{bad}$(vi)$ and the fact that $1+\la + \la^2 + \ldots = b_2$. The same argument as in Case \ref{c2c3=1} allows us to assume that $j \leq b_2$, and this gives an immediate upper bound of $\frac{3}{2b_2} + \frac{1}{4} + \frac{\la^{j-5}}{4}< \frac{3}{2b_2}+\frac{1}{2}$ for \rrr{jw2}. This is clearly bounded by $1$ for $b_2 \geq 3$.

\vski

Now suppose $12 < D \leq 16$. By Lemma 5.2 of \cite{hbk}, $b_3 \geq 2 c_4$, so that $\frac{b_3}{c_3} \geq 2$. We also have $j \geq 5$ by Lemma \ref{bad}$(iii)$. We have, setting $\la = \frac{b_2-1}{b_2}$ as before,

\bea \lll{jw3}
&& \frac{\phi_2+ \ldots + \phi_{15}}{\phi_1}
\\ \nn && \hspace{1cm} < \frac{1}{b_2} + \frac{1}{2b_2} + \frac{\la}{2b_2}+ \frac{\la^2}{2b_2} +\ldots + \frac{\la^{j-4}}{2b_2} + \phi_{j} + \ldots + \phi_{15}
\\ \nn && \hspace{1cm} = \frac{1}{b_2} + \frac{1}{2}(1-\la^{j-3}) + \phi_{j} + \ldots + \phi_{15}
\eea

We now claim $\phi_{j} + \ldots + \phi_{15} < 7 \phi_{j-1}$. This is clear if $j \geq 9$, since $\phi_{j} + \ldots + \phi_{15}$ contains at most $7$ elements, whereas Lemma \ref{bad}$(vi)$ applies if $j \leq 8$. Since $\phi_{j-1} < \frac{\la^{j-4}}{2b_2}$, \rrr{jw3} gives

\be \lll{jw4}
\frac{\phi_2+ \ldots + \phi_{15}}{\phi_1} < \frac{7 \la^{j-4}+2-b_2\la^{j-3}}{2b_2} + \frac{1}{2} = \frac{(8-b_2) \la^{j-4}+2}{2b_2} + \frac{1}{2}
\ee

where the identity $b_2 \la = b_2-1$ was used. It is now immediate that this is bounded by $1$ if $b_2 \geq 5$, while we also obtain $1$ as a bound if we substitute $b_2=4, j \geq 7$. Suppose $b_2=4, j=6$. Then by Lemma \ref{bad}$(vi)$, $\phi_{6} + \ldots + \phi_{15} < 5 \phi_{5}$, and proceeding as before we obtain

\be \lll{jw5}
\frac{\phi_2+ \ldots + \phi_{15}}{\phi_1} < \frac{5 \la^{2}+2-b_2\la^{3}}{2b_2} + \frac{1}{2} = \frac{(6-b_2) \la^{2}+2}{2b_2} + \frac{1}{2} ,
\ee

and this is less than $1$. Now assume $b_2=4,j=5$. We apply Lemma \ref{bad}$(vi)$ once again to get

\bea \lll{jw6}
&& \frac{\phi_2+ \ldots + \phi_{15}}{\phi_1}
\\ \nn && \hspace{1cm} < \frac{1}{b_2} + \frac{1}{2b_2} + \frac{\la}{2b_2} + \phi_{j} + \ldots + \phi_{13}
\\ \nn && \hspace{1cm} < \frac{3}{2b_2} + 4 \times \frac{\la}{2b_2} = \frac{7}{2b_2}
\eea

This is clearly bounded by $1$, since $b_2 \geq 4$. This completes the proof for $12<D \leq 16$.

\vski

Now suppose $D\leq 12$. Case \ref{niceD} covers $b_2 \geq 5$, so we may assume $b_2=4$. Case \ref{niceD} also covers $b_2=4, D \leq 10$, so we may assume $D \geq 11$, whence $j \geq 5$ by Lemma \ref{bad}$(iii)$. If $c_3 \geq 3$, then by Lemma 5.2 of \cite{hbk} we would have $D \leq 8$. Thus, $c_3=c_3=2$, and $b_2$ may be $3$ or $4$. Let us suppose first that $j=5$. There are two possibilities to consider. $b_3$ may be $3$, in which case $b_4=3, c_4=2$ as well and we may use \rrr{comm} and Lemma \ref{bad}$(vi)$ to obtain

\be \lll{jw7}
\frac{\phi_2+ \ldots + \phi_{11}}{\phi_1}< \frac{1}{4} + \frac{2}{4 \times 3} + 5 \times \frac{2^2}{4 \times 3^2} < 1
\ee

If that does not occur, then $b_3 = 4$, and $\frac{c_4}{b_4} \leq \frac{3}{4}$, so we have

\be \lll{jw8}
\frac{\phi_2+ \ldots + \phi_{11}}{\phi_1}< \frac{1}{4} + \frac{2}{4 \times 4} + 5 \times \frac{2\times 3}{4 \times 4 \times 3} = 1
\ee

$j=5$ is thereby handled. Assume now that $j \geq 6$. Note that if $b_3=3$ then since it was shown earlier that $c_5 \geq 3$ we have $c_5 \geq b_5$ and thus $j=5$. We may therefore assume $b_3=4$. Note that $\frac{c_i}{b_i} \leq \frac{3}{4}$ for $i < j$, so we have

\bea \lll{jw9}
&& \nn \frac{\phi_2+ \ldots +\phi_{11}}{\phi_1}
\\ \nn && \hspace{1cm} < \frac{1}{4} + \frac{1}{8} + \frac{1}{8} \times \frac{3}{4} +  \frac{1}{8} \times \frac{3^2}{4^2} + \ldots + \frac{1}{8} \times \frac{3^{j-4}}{4^{j-4}} + \frac{\phi_{j} + \ldots + \phi_{11}}{\phi_1}
\\ \nn && \hspace{1cm} < \frac{1}{4} + \frac{1}{8} + \frac{1}{8} \times \frac{3}{4} +  \frac{1}{8} \times \frac{3^2}{4^2} + \ldots + \frac{1}{8} \times \frac{3^{j-4}}{4^{j-4}} + (12-j) \times \Big(\frac{1}{8} \times \frac{3^{j-4}}{4^{j-4}}\Big)
\\ \nn && \hspace{1cm} = \frac{1}{4} + \frac{1}{2}(1-\frac{3^{j-3}}{4^{j-3}}) + (12-j) \times \Big(\frac{1}{8} \times \frac{3^{j-4}}{4^{j-4}}\Big)
\\ \nn && \hspace{1cm} = \frac{3}{4} + \frac{(9-j)(3/4)^{j-4}}{8}
\eea

The last expression is clearly either negative or decreasing in $j$, and plugging in $j=6$ gives a bound less than $1$. This completes the proof for $b_2=4$, and completes Case \ref{c2=1}. \hfill $\triangle$

\vski

\ccases{b_2<=3}: $b_2\leq3$, $c_2=1$

\vski

Suppose first that $a_1=a_2$. As $c_2=c_1=1$, we have $b_1=b_2$. Also, we know that $a_1+1$ divides $b_2$, as $a_1+1$ divides $k=1+a_1+b_1$ and $b_1=b_2$. i.e., $\alpha\times(a_1+1)=b_2\leq3$ for some integer $\alpha$. If $b_1=b_2\in\{1,2\}$, then $a_1=a_2\in\{0,1\}$, so the valency $k$ is at most $4$, and we may appeal to Case \ref{deg3} or Case \ref{deg4}. We may therefore assume $b_1=b_2=3$. Then we know $a_1=a_2\in\{0,2\}$. But $a_1=a_2=0$ implies that $k=4$, handled in Case \ref{deg4}, so we may assume $a_1=a_2=2$. This implies that $k=6$ and $G$ is a line graph, as $c_2=1$ and $a_1=\frac{1}{2}k-1$. Then the graph is either the flag graph of a regular generalized $D$-gon of order $(3,3)$ or the line graph of a Moore graph by \cite[Theorem 4.2.16]{drgraphs}. It is straightforward to check that these graphs satisfy the desired inequality.

\vski

Now let us consider the case $a_1\neq a_2$, which implies $a_2>a_1$. First, we will assume $c_3=1$. Let $x$ and $y$ be vertices of the graph $G$, and put $\Delta(x,y):=\{x\}\cup A\cup C$, where $C$ is the component of $K_2(x)$(the set of vertices at distance two from $x$) containing $x$, and $A$ is the union of all maximal cliques containing $x$ and a neighbor of a point of $C$. Then $\Delta(x,y)$ is a strongly regular graph with parameters $(v',k',a_1',c_2')=(v',a_2+1,a_1,1)$ by \cite[Proposition 4.3.11]{drgraphs}. As the size of a maximal clique in $G$ is at most $b_2+c_3$, we have $a_1+2\leq b_2+c_3$, which implies $a_1+1\leq b_2\leq3$. By \cite[Corollary 4.3.12]{drgraphs}, we also have $b_2\geq a_2-1$. Then $k'=a_2+1\leq a_1+b_2+1\leq2b_2\leq6$ holds. As $k'\leq6$ and $b_2\leq3$, we know $v'\leq25$. By checking the table of strongly regular graphs (http://www.win.tue.nl/$\sim$aeb/graphs/srg/), we find that the only possible case is that the strongly regular graph $\Delta(x,y)$ is the pentagon. i.e., $G$ has valency $k=5$ and intersection numbers $a_1=0$, $a_2=1$ and $b_2=3$. There is no distance-regular graph with $c_4=1$ by \cite[Theorem 1.1]{hiraki} and \cite[Theorem 1.2]{chk}, so we may assume $c_4 \geq 2$. If $c_3 \geq 2$, then by Lemma 5.2 of \cite{hbk} we have $D \leq 8$, and we may then appeal to Case \ref{niceD}. We may therefore assume $c_3=1$, so we have $c_2=c_3=1, c_4 \geq 2$. Assume first that $b_3=3$ as well. Then $\frac{c_i}{b_i} \leq \frac{2}{3}$ for $i < j$, so using Lemma \ref{bad}$(vi)$ we write

\bea \lll{jc}
&& \nn \frac{\phi_2+ \ldots \phi_{D-1}}{\phi_1}
\\ \nn && \hspace{1cm} < \frac{1}{3} + \frac{1}{9} + \frac{1}{9} \times \frac{2}{3} +  \frac{1}{9} \times \frac{2^2}{3^2} + \ldots + \frac{1}{9} \times \frac{2^{j-4}}{3^{j-4}} + \frac{\phi_{j} + \ldots + \phi_{D-1}}{\phi_1}
\\ \nn && \hspace{1cm} < \frac{1}{3} + \frac{1}{9} + \frac{1}{9} \times \frac{2}{3} +  \frac{1}{9} \times \frac{2^2}{3^2} + \ldots + \frac{1}{9} \times \frac{2^{j-4}}{3^{j-4}} + (j-1/2) \times \Big(\frac{1}{9} \times \frac{2^{j-4}}{3^{j-4}}\Big)
\\ \nn && \hspace{1cm} = \frac{1}{3} + \frac{1}{3}(1-\frac{2^{j-3}}{3^{j-3}}) + (j-1/2) \times \Big(\frac{1}{9} \times \frac{2^{j-4}}{3^{j-4}}\Big)
\\ \nn && \hspace{1cm} = \frac{2}{3} + \frac{(j-5/2)(2/3)^{j-4}}{9}
\eea

It is simple to verify as in Case \ref{c2c3=1} that the last expression is decreasing in $j$ for $j \geq 4$, and plugging in $j=4$ gives a bound of $\frac{5}{6}$. Now let us assume $b_3\leq2$. As $c_4\geq2$, we know $c_7>c_4\geq2\leq b_3$. So, the diameter $D$ is at most nine and the number of vertices, $n$, is at most $446$. Consulting the table of intersection arrays of \cite[Chapter 14]{drgraphs}, we see that there are no distance-regular graph with $k=5$, $b_1=4$, $b_2=3$, $b_3\leq2$, $c_2=c_3=1$, $c_4\geq2$, $D\leq9$ and $v\leq446$.  \hfill $\triangle$

\vski

The reader may now verify that we have covered all cases, and the proof of Theorem \ref{thispap} is complete. \qed

\section{The strongly regular case} \lll{srgraphs}

Let $G$ be a strongly regular graph with diameter $D=2$(if $D=1$ then $G$ is a complete graph, which is a trivial case for our considerations). Suppose $G$ has parameters $(v,k,a_1,c_2)$. Here note that if $c_2=k$ then $G$ is complete multipartite (or, see \cite[Theorem 1.3.1 $(v)$]{drgraphs}).

\begin{theorem} \lll{sreg}
Let $G$ be a strongly regular graph with diameter $2$ with parameters $(v,k,a_1,c_2)$, and assume that $G$ is not complete multipartite. Then $b_1=k-a_1-1\geq {\rm min}\{\frac{5}{16}k, \frac{2}{1+\sqrt{2}}\sqrt{k}\}$. In particular, there are only finitely many co-connected, non-complete strongly regular graphs with fixed $b_1$.
\end{theorem}

{\bf Proof:} Let $\theta_1>\theta_2$ be non-trivial eigenvalues of $G$. If one of $\theta_1$ and $\theta_2$ is non-integral, then both are non-integral and they have the same multiplicity. i.e., $G$ is a conference graph and hence $G$ has the parameters $(4b_1+1,2b_1,b_1-1,b_1)$. Thus, $k=2b_1$ in this case, and the theorem holds. From now on, we may assume that both $\theta_1$ and $\theta_2$ are integral. Note that we may assume that $\theta_1\geq1$ and $\theta_2\leq-3$, unless $G$ has $b_1\geq\frac{5}{16}k$ by \cite[Theorem 1.3.1 and Theorem 3.12.4]{drgraphs}. By \cite[Theorem 3.6]{kpy}, we have $-b_1=(\theta_1+1)(\theta_2+1)$. Then we know that $\theta_1\leq \frac{b_1}{2}-1$ and $\theta_2\geq-\frac{b_1}{2}-1$ hold. As $G$ is not complete multipartite it is not antipodal, hence $(k_2)^2>k_2(k_2-1)\geq k$ holds by \cite[Proposition 5.6.1]{drgraphs}. i.e., $\sqrt{k}b_1>c_2$ holds. As $c_2=\theta_1\theta_2+k$ (see \cite[Chapter 10]{gr}), we have $-\sqrt{k}b_1<-\theta_1\theta_2-k<\frac{b_1^2}{4}-k$. The quadratic inequality in $\sqrt{k}$, $k-b_1\sqrt{k}-\frac{b_1^2}{4}<0$, gives the inequality  $\sqrt{k}<\frac{1+\sqrt{2}}{2}b_1$. \qed

A complete multipartite graph $K_{m \times t}$ can be realized in the following way. Let the vertices be the set of all pairs $(a,b)$, with $a \in \{1, 2, \ldots , m\}$ and $b \in \{1, 2, \ldots , t\}$. We let $(a,b) \sim (a',b')$ if and only if $a \neq a'$. This is a strongly regular graph containing $mt$ vertices, and $k=(m-1)t, b_1=t-1$. Theorem \ref{sreg} therefore fails for this graph, since $m$ may be arbitrarily large. However one may readily check that $\phi_0=mt-1, \phi_1 = 1$, so $\frac{r_2}{r_1} = \frac{\phi_0+\phi_1}{\phi_0} = 1+\frac{1}{mt-1} \leq 1+\frac{1}{k}$. The following is therefore a corollary of the preceding theorem, together with the fact that $\phi_1 < \frac{\phi_0}{b_1}$.

\begin{corollary} \lll{last}
Let $C(k)={\rm min}\{\frac{5}{16}k, \frac{2}{1+\sqrt{2}}\sqrt{k}\}$. Then, if $G$ is a strongly regular graph with diameter $2$, we have

\be \label{}
\frac{r_2}{r_1} = \frac{\phi_0+\phi_1}{\phi_0} < 1 + \frac{1}{C(k)}
\ee
\end{corollary}

This, combined with Theorem \ref{latcor}, validates the statement given in the abstract and introduction, that all points are more or less equidistant under the electric resistance metric in distance-regular graphs with large valency.

\section{Concluding remarks} \lll{conc}

We would like to remark first that Theorem \ref{bigguy}, the subject of \cite{markool}, is almost a corollary of Theorem \ref{thispap}. This is because if $b_1 \geq 2$ then applying Theorem \ref{thispap} gives

\be \label{}
\frac{\phi_0+ \ldots \phi_{D-1}}{\phi_0} \leq \frac{\phi_0 + 2\phi_1}{\phi_0} < \frac{\phi_0 + 2\times \frac{\phi_0}{b_1}}{\phi_0} \leq 2
\ee

Improving the constant slightly from $2$ to $K=1+\frac{94}{101}$ and dealing with the case $b_1=1$ is not difficult. Indeed, the proof of Theorem \ref{bigguy} given in \cite{markool} can be made substantially shorter and cleaner using some of the new ideas in this paper, in particular Lemma \ref{bad}$(vi)$.

\vski

Next, we would like to give examples of some intersection arrays which are ruled out by Theorem \ref{thispap}.

\vski

\noi \begin{tabular}{ l c r }
Intersection array & Vertices & $\frac{\phi_2+ \ldots \phi_{D-1}}{\phi_1}$ \\
\hline
(3, 2, 2, 1, 1, 1, 1;1, 1, 1, 1, 1, 1, 3) & 62 & 1.2069 \\
(5, 2, 2, 1, 1, 1, 1;1, 1, 1, 1, 1, 1, 4) & 101 & 1.18421 \\
(6, 4, 4, 3, 3, 2, 1, 1;1, 2, 2, 3, 3, 3, 4, 4) & 112 &      1.03333\\
\end{tabular}

\vski

We hope that Theorem \ref{thispap} might be found useful for ruling out possible intersection arrays, as well as for proving other facts.

\section{Acknowledgements}

The first author was partially supported by the Basic Science Research Program through the National Research
Foundation of Korea(NRF) funded by the Ministry of Education, Science and Technology (Grant \# 2009-0089826). The second author was supported by Priority Research Centers Program through the National Research Foundation of Korea (NRF) funded by the Ministry of Education, Science and Technology (Grant \#2009-0094070).


\begin{thebibliography}{1}

\bibitem{hbk}
S. Bang, A. Hiraki, J. Koolen (2006) {\it Improving diameter bounds for distance-regular graphs}, Eur. J. Comb. 27(1), p. 79-89

\bibitem{biggs2}
N. Biggs (1993) {\it Potential theory on distance-regular graphs.} Combinatorics, Probability and Computing 2, p. 243-255.

\bibitem{biggs}
N. Biggs (1997) {\it Algebraic Potential Theory on Graphs}, Bulletin of the London Mathematical Society, 29(6), p. 641-682.

\bibitem{drgraphs}
A. Brouwer, A. Cohen, A. Neumaier (1989) {\it Distance Regular Graphs}, Springer-Verlag.

\bibitem{broukool}
A. Brouwer, J. Koolen (1999) {\it The Distance-Regular Graphs of Valency Four}, Journal of Algebraic Combinatorics, 10(1), p. 5-24.

\bibitem{chk}
Y.L. Chen, A. Hiraki, J.H. Koolen (1998) {\it On distance-regular graphs with $c_4=1$ and $a_1\neq a_2$}, Kyushu J. Math. 52, p. 265-277.

\bibitem{gr}
C. Godsil, G. Royle (2001) {\it Algebraic Graph Theory}, Springer-Verlag, Berlin.

\bibitem{hiraki}
A. Hiraki, (1995) {\it Distance-regular subgraphs in a distance-regular graph III}, Europ. J. Combin. 17, p. 629-636.

\bibitem{koolen}
J. Koolen, (1992) {\it On subgraphs in distance-regular graphs}, J. Algebraic Combin. 1(4), p. 353-362.

\bibitem{parkool2}
J. Koolen, J. Park (2010) {\it Distance-regular graphs with large a1 or c2}, arXiv:1008.1209v1.

\bibitem{kpy}
J. Koolen, J. Park, H. Yu (2010) {\it An inequality involving the second largest and smallest eigenvalue of a distance-regular graph}, arXiv:1004.1056v1.

\bibitem{markool}
G. Markowsky, J. Koolen (2010) {\it A Conjecture of Biggs Concerning the Resistance of a Distance-Regular Graph}, Electronic Journal of Combinatorics, v. 17(1).

\bibitem{parkool}
J. Park, J. Koolen, G. Markowsky (2010) {\it There are only finitely many distance-regular graphs with valency $k$ at least three, fixed ratio $\frac{k_2}{k}$ and large diameter}, preprint, arXiv:1012.2632v1.







\end{thebibliography}
\end{document}